\documentclass[11pt, oneside]{article}   	
\usepackage{geometry}                		
\usepackage[parfill]{parskip}    			
\usepackage{cancel}
\usepackage{graphicx}				
\usepackage{amssymb}

\begin{document}

\title{Patterns Within the Markov Tree}

\author{Robert A. Gore\\Cheyney, PA}

\maketitle

\begin{abstract}
An analysis of the Markov tree is presented. Markov triplets, $\{x, R, z\}$, are the positive integer solutions to the Diophantine equation $x^2+R^2+z^2=3x R z.$ Inspired by patterns of the Fibonacci and Pell triplets in \textsc{Region 1} and \textsc{Region 2} of the tree, an investigation of interior regions of the Markov tree finds generating functions and sequence functions for all triplets of all regions. These sequence functions lead to the discovery of a Pell equation for the Markov region numbers along the edges of all regions. Analysis of this Pell equation leads to the resolution of the Uniqueness Conjecture. Further analysis using these sequence functions finds palindromic repeat cycles of the last digits of region numbers along the edges of all regions. 
  
  Then, since all Markov numbers are the sum of the squares of two integers and again inspired by the patterns of the two unique squares which sum to form the region numbers of certain Fibonacci triplets in \textsc{Region 1}, an investigation of interior regions of the Markov tree finds generating functions and sequence functions for the two special square terms which sum to form the region numbers of the triplets along the edges of all regions. Further analysis using these sequence functions finds palindromic repeat cycles of the last digits of these two special square terms for all regions. 
\end{abstract}

\section{The Interconnectivity of the Markov Tree}
Markov triplets $\{x, R, z\}$ are integer solutions to the Diophantine equation: 
$$x^2+R^2+z^2=3x R z.$$
Define the first three solutions, $\{x, R, z\}$, as the \textit{ordered} sets $\{1, 1, 1\}$, $\{1, 2, 1\}$ and $\{1, 5, 2\}$ and define $R$ with $R>x$ and $R>z$ as the region number of triplet $\{x, R, z\}$. The first two of these solutions are singular triplets since their members are not unique. The last is non-singular. It is well known that every non-singular Markov triplet is connected to every other non-singular Markov triplet, since any parent triplet $\{x, R, z\}$  is related to an \textit{ordered} left hand child triplet $\{x, 3Rx-z, R\}$ and an \textit{ordered} right hand child triplet $\{R, 3Rz-x, z\}$. This ordered notation for the children of a parent allows the parent triplet to be deduced from a child triplet. The parent of a triplet $\{x, R, z\}$ with $x<z$ is $\{x, z, 3 x z-R\}$ and the parent of a triplet $\{x, R, z\}$ with $x>z$ is $\{3 x z-R,x,z\}$. 

\subsection{Fibonacci and Pell Regions: Recurrence Kernels}
Consider \textsc{Region~1} and \textsc{Region~2} of the Markov tree. The region numbers of the first few triplets of the Markov tree are shown in Fig.\ 1.
\begin{figure}[ht]
\begin{center}
\includegraphics[width=1.0\columnwidth]{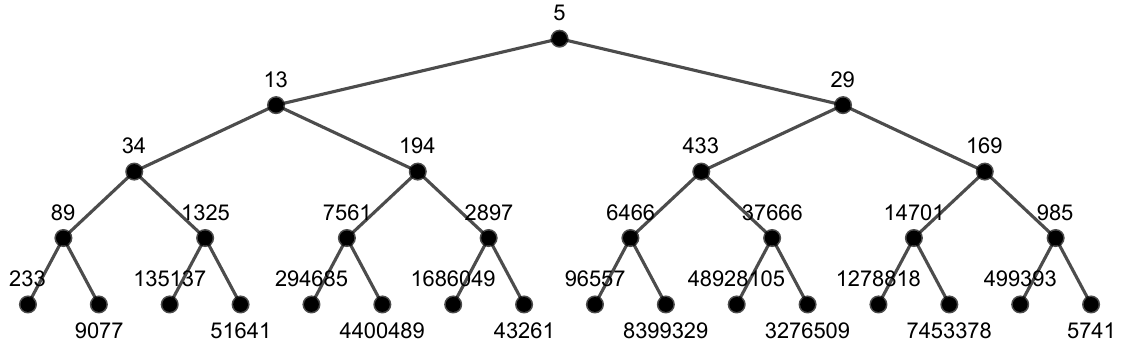}
\end{center}
\caption{The Markov Tree Starting with the Triplet $\{1, 5, 2\}$. }
\end{figure}
Although the triplets along the left and right outer edges of Fig.\ 1 are related to the odd-indexed Fibonacci and Pell numbers, they have a deeper symmetry more closely related to the Markov tree. The region numbers of \textsc{Region~1} of the Markov tree, \(\{13, 34, 89, ...\}\), have a recurrence kernel of \{3*1,~-1\}, while the region numbers of \textsc{Region~2} of the Markov tree, \(\{29, 169, 985, ...\}\), have a  recurrence kernel of \{3*2,~-1\}. Note that the first terms of these  recurrence kernels are in the 1:2 ratio of their respective \textsc{Regions} \{1, 2\}. Since the single string of Fibonacci related triplets are all right hand Markov children and the single string of Pell related triplets are all left hand Markov children, it is useful to consider left hand children \textit{separately} from right hand children when dealing with the two strings of \textit{interior} region triplets of the Markov tree. 

For example, consider the first internal region of the Markov tree, \textsc{Region~5}. This region is headed by the triplet \(\{1, 5, 2\}\). Both the region numbers along the left edge, \(\{13, 194, 2897, ...\}\), and the region numbers along the right edge, \(\{29, 433, 6466, ...\}\), have a  recurrence kernel of \{3*5,~-1\}. 
Note that the pattern of  recurrence kernels for these three \textsc{Regions} \{1, 2, 5\} is: 
$$\{\{3*1,~-1\},\{3*2,~-1\},\{3*5,~-1\}\}.$$
The next section extends these results for \textit{all} regions. It starts by deriving the recurrence kernel for the region numbers along the left and right edges of an arbitrary interior Markov region headed by the triplet $\{x, R, z\}$. Next, generating functions and sequence functions are derived for these left and right edge triplets.

\subsection{Generating Functions and Sequence Functions for All Triplets}
Define $\{x, R, z\}$ with $R>x$ and $R>z$ as the triplet at the head of an arbitrary interior \textsc{Region~$R$} $(R\ge5)$ of the Markov tree. Evaluate the first three left and right edge children of this triplet and extract their region numbers:\\
$\{\textsc{lhs}_1, \textsc{lhs}_2, \textsc{lhs}_3\} = \{3Rx-z, ~(-1 + 9 R^2) x - 3 R z, ~(1-9R^2)z-3R x(2-9R^2)\}$\\
$\{\textsc{rhs}_1, \textsc{rhs}_2, \textsc{rhs}_3\}=\{3Rz-x, ~(-1 + 9 R^2) z - 3 R x, ~(1-9R^2)x-3R z(2-9R^2)\}.$\\
The  recurrence kernel of both of these sets of region numbers is $\{3R,~-1\}.$
Since no assumptions were made with respect to the interior triplet, $\{x, R, z\}$, and only the formulas for the Markov children were used,  $\{3R,~-1\}$ is the  recurrence kernel for the region numbers of left and right edge triplets of any interior Markov Region~$R$. The intermediate triplet members, Max\{x,z\}, have this same recurrence kernel since they are the region numbers of their parents.
Thus, the ordered generating functions of the left and right edge triplets of region \(\{x, R, z\}\) are:
$$
\textsc{H}_{gf} [\{x, R, z\},n] = 
\left(
\begin{array}{cc}
 \{\frac{x-n z}{1-3R n+n^2},   ~\frac{(3R x-z)-n x}{1-3R n+n^2},   ~\frac{R}{1-n}\}  & \textsc{left edge } \\
 \\
\{\frac{R}{1-n},    ~\frac{(3R z-x)-n z}{1-3R n+n^2},   ~\frac{z-n x}{1-3R n+n^2}\} & \textsc{right edge}\\
\end{array}
\right.
$$
It is useful to have sequence functions for the $n^{th}$ member of these triplets. For $n \ge 1$, the ordered sequence functions of the left and right edge triplets of region \(\{x, R, z\}\) are:
$$
\textsc{H}_{sf} [\{x, R, z\},n] = 
\left(
\begin{array}{ll}
\{U[\{z, R, x\},n] , ~U[\{z, R, x\},n+1] , ~R\}   & \textsc{left edge}\\
\{R, ~U[\{x, R, z\},n+1] , ~U[\{x, R, z\},n]\}    & \textsc{right edge}\\
\end{array}
\right.
$$
Here, $U[\{a, R, b\},n] $ is a Lucas sequence defined as:
$$U[\{a, R, b\},n]    = b~U_{n}(3R,1) - a~U_{n-1}(3R,1)$$
with initial values $$U[\{a, R, b\},0] = a$$ and $$U[\{a, R, b\},1] = b.$$

The sequence $U_{k}(3R,1)$ and its associated sequence $V_{k}(3R,1)$ are specific cases of the general Lucas sequences $U_{k}(P, Q )$ and $V_{k}(P, Q )$ given by:
$$U_{k}(3R,1) = \frac{\sinh(k \theta)} {\sinh(\theta)}$$
$$V_{k}(3R,1) = 2 \cosh(k \theta)$$
with $ $cosh$ (\theta) = \frac{3R}{2}.$

The sequence function, $\textsc{H}_{sf} [\{x, R, z\},n]$, has several interesting properties. 
\begin{itemize}
\item{$\textsc{H}_{sf} [\{x, R, z\},n]$ \textbf{Works for all Regions}

Even though these expressions were derived for an arbitrary interior region, they solve the Markov equation for all triplets of all regions. For \textsc{Region~$1$}, headed by $\{1, 1, 1\}$, the right edge formulas are required as these Fibonacci triplets are right hand children. Similarly, for \textsc{Region~$2$}, headed by $\{1, 2, 1\}$, the left edge formulas are required as these Pell triplets are left hand children.}

\item{$\textsc{H}_{sf} [\{x, R, z\},n]$ \textbf{Works for All Integer $n$}

Even though these expressions were derived for positive $n$, they have well defined values for all $n$:
 $U_n(3 R,1)$ is an odd function of $n$ while $V_n(3 R,1)$ is even.
\mbox{Table 1} shows the output of the \{left, right\} edge triplet region numbers for positive and negative $n$. 
 Note that the \{left, right\} edge expression for negative $n$ produces the values along the opposite edge for positive $n$. Each edge expression has the full information content of all triplets along both edges although only two triplets along one edge were sufficient to derive the sequence functions. This behavior occurs for all interior regions. 
\begin{table}[ht]
\begin{center}
\caption{Triplet Region Numbers Derived from $\textsc{H}_{sf} [\{x, R, z\},n]$}
\begin{tabular}{|c|c|c|c|}
\hline
Index	&Left Edge Triplet Region 		&Right Edge Triplet Region\\
$n$		&Numbers Derived from		&Numbers Derived from\\
		&$\textsc{H}_{sf} [\{x, R, z\},n]$	&$\textsc{H}_{sf} [\{x, R, z\},n]$\\
\hline
 -4 & $z U_4(3 R,1)-x U_3(3 R,1)$ & $x U_4(3 R,1)-z U_3(3 R,1)$ \\
 -3 & $z U_3(3 R,1)-x U_2(3 R,1)$ & $x U_3(3 R,1)-z U_2(3 R,1)$ \\
 -2 & $z U_2(3 R,1)-x$ & $x U_2(3 R,1)-z$ \\
 -1 & $z$ & $x$ \\
 0 & $x$ & $z$ \\
 1 & $x U_2(3 R,1)-z$ & $z U_2(3 R,1)-x$ \\
 2 & $x U_3(3 R,1)-z U_2(3 R,1)$ & $z U_3(3 R,1)-x U_2(3 R,1)$ \\
 3 & $x U_4(3 R,1)-z U_3(3 R,1)$ & $z U_4(3 R,1)-x U_3(3 R,1)$ \\
 4 & $x U_5(3 R,1)-z U_4(3 R,1)$ & $z U_5(3 R,1)-x U_4(3 R,1)$ \\
\hline
\end{tabular}
\end{center}
\end{table}
}

\item{\textbf{The Grandparent of} $\textsc{H}_{sf} [\{x, R, z\},n]$ \textbf{is the Secondary Solution of the Markov Equation for Fixed $x$ and $z$}}

Transform the sequence function, $\textsc{H}_{sf} [\{x, R, z\},n],$ for a right edge triplet,
$$\{R, ~U[\{x,R,z\},n+1] , ~U[\{x,R,z\},n]\},$$
by swapping the two terms which involve a Lucas sequence and shifting $n \rightarrow (n-1)$ to obtain the secondary solution of the Markov equation for fixed $x$ and $z$,
$$\{R, ~U[\{x,R,z\},n-1], ~U[\{x,R,z\},n]\}.$$
This works for all non-singular triplets. 
In what follows, the secondary solution of the Markov equation for fixed $x$ and $z$ is termed the region sibling number.
Additionally, triplets hanging outside off of a region $R$ have a sibling number of $R$.

\end{itemize}

\section{A Proof of  The Uniqueness Conjecture}
This analysis is patterned after the standard analysis of the Fibonacci sequence. It is well known that the Fibonacci numbers, $F_n = U_n(1,-1)$, and the associated Lucas numbers, $L_n = V_n(1,-1)$, obey the Pell equation
$$L_n^{2} - 5 ~F_n^{2} =4(-1)^n.$$
Using standard techniques for solving a Pell equation of the form $\Phi_n^2 -5 \Psi_n^2 = 4(-1)^n,$
this Pell equation has been shown to have solutions \textit{if and only if} $\Phi_n$  is the Lucas number $L_n$ and $\Psi_n$  is the Fibonacci number $F_n.$
Thus, an arbitrary integer \textbf{J} is a Fibonacci number \textit{if and only if} either or both $\{\sqrt{5 \textbf{J}^2 + 4},\sqrt{5 \textbf{J}^2 - 4}\}$ is a perfect square. More specifically, an arbitrary integer \textbf{J} is an \textit{even-indexed} Fibonacci number \textit{if and only if} $\{\sqrt{5 \textbf{J}^2 + 4}$ is a perfect square. And finally, an arbitrary integer \textbf{J} is the $n^{th}$ Fibonacci number \textit{if and only if} $L_n^{2} - 5 \textbf{J}^2 =4(-1)^n.$ 

This analysis can be generalized to form a proof of the Fr\"obenius Uniqueness Conjecture as follows:
Recall the sequence function for the region number of the right edge triplet of an arbitrary Markov region $\{x, R, z\},$
$$U[\{x,R,z\},n+1] = z U_{n+1}(3R,1) - x U_{n}(3R,1),$$
and define its associated Lucas sequence by
$$V[\{x,R,z\},n+1] = z V_{n+1}(3R,1) - x V_{n}(3R,1).$$
These two sequences satisfy a particular Pell equation
$$ V[\{x,R,z\},n+1] ^2 - D(R) ~U[\{x,R,z\},n+1] ^2 = - (2R)^2$$
for all Markov triplets and all integer $n$ with the definition of the discriminant
$$D(R)=  (3 R)^2 - 4.$$
To verify this, expand $V[\{x,R,z\},n+1]$ and $U[\{x,R,z\},n+1]$, recalling that 
$$U_{k}(3R,1) = \frac{\sinh(k \theta)} {\sinh(\theta)}$$
$$V_{k}(3R,1) = 2 \cosh(k \theta).$$
And since $\cosh (\theta) = \frac{3R}{2},$ replace $D(R)$ with $4\sinh^2(\theta).$ 
Collecting $x$, $R$, and $z$ terms, this Pell equation becomes the Markov equation after replacing 
the $\cosh(\theta)$ term:
$$4 x^2 [\cosh^2(n \theta) - \sinh^2(n \theta)] + 4 R^2 +4 z^2 [\cosh^2((n+1) \theta) - \sinh^2((n+1) \theta)] = 8 x z \cosh(\theta).$$
Now, one can prove that $V[\{x,R,z\},n+1]$ and $U[\{x,R,z\},n+1]$ are the \textit{only} solutions to this Pell equation by choosing a convenient baseline solution to 
$$\textbf{X}_{n}^2 - D(R) \textbf{Y}_{n}^2=- (2R)^2$$
as
$$\{\textbf{X}_0, \textbf{Y}_0\} = \{V[\{x,R,z\},1], U[\{x,R,z\},1]\}$$
and by noting that the simplest solution to the unit residual Pell equation
$$k_{n}^2 - D(R)~ l_{n}^2=1$$
is
$$\{k_0, l_0\}=\frac{1}{2}\{3R,1\}.$$
Thus, the set of solutions to the Pell equation,
$$\textbf{X}_{m}^2 - D(R) \textbf{Y}_{m}^2=- (2R)^2,$$ (referenced to the baseline solution) is given by
$$(\textbf{X}_m + \sqrt{D(R)} ~\textbf{Y}_m) = (k_0+\sqrt{D(R)} ~l_0)^m ~(\textbf{X}_0+\sqrt{D[R]}~\textbf{Y}_0).$$
Since the results, $\{\textbf{X}_m, \textbf{Y}_m\}$, have the same recurrence kernels, generating functions, and values as $\{V[\{x,R,z\},m+1], ~U[\{x,R,z\},m+1]\}$ respectively, these are the \textit{only} solutions:
$$\textbf{K} ^2 - D(R) ~\textbf{J} ^2 = - (2R)^2$$ is true  \textit{if and only if} \textbf{K} = $V[\{x,R,z\},n+1]$ and \textbf{J} = $U[\{x,R,z\},n+1]$. 

For clarity, this analysis was designed to treat left and right edges separately. Yet, the above right edge analysis can be transformed to apply to left edges by simply substituting the word ``left" for ``right"  and interchanging $x \leftrightarrow z$.

 A better way to extend the above right edge analysis to apply to left edges is to  recall that Table 1 of Section 1.2 showed that the general expression $U[\{x,R,z\},n]$ wraps around and is valid for \textit{all} integer $n$. Since nothing within the above right edge analysis restricted the value of $n$, it applies for \textit{all} $n$.
Thus, the above right edge analysis, which uses $n=\{1, 2, 3...\}$ for the right edges, can also be directly used with $n=\{-2,-3,-4...\}$ for the left edges as shown in Table 1 of Section 1.2.
Notice that two values of $n$, $\{n=0, n=-1\}$, are missing here but not in Table 1 of Section 1.2 \textit{where the values for these missing cases are $x$ and $z$.}
\newline
To highlight the significance of these missing cases, imagine using software which can solve integer domain equations such as the one parameter Pell equation,
$$\textbf{K}^2 - [(3R)^2-4] ~\textbf{J}^2= - (2R)^2,$$
for different values of the parameter $R$. Such an analysis using \textit{Mathematica} 
shows that if and only if the parameter $R$ is equal to a Markov number does a solution exist and that the \textit{two smallest solutions, \textbf{J}, for a given parameter $R$  are always the $x$ and $z$ values which complete the Markov triplet $\{x, R, z\}$ for region number $R$.} 

Table 2 shows the results of tabulating the first few \textit{Mathematica} solutions to this Pell equation for those values of the parameter $R$ which have a solution.
 
 \begin{table}[ht]
\begin{center}
\caption{Solutions, \textbf{J}, of a Pell Equation for the Parameter $R$}
\begin{tabular}{|c|c|}
\hline
Parameter $R$					&Sorted Positive Integer Solutions\\
\hline
\hline
 1 & \{1, 2, 5, 13, 34, 89 ...\} \\
 2 & \{1, 5, 29, 169, 985, 5741 ...\} \\
 5 & \{1, 2, 13, 29, 194, 433 ...\} \\
 13 & \{1, 5, 34, 194, 1325, 7561 ...\}\\
 29 & \{2, 5, 169, 433, 14701, 37666 ...\}\\
 34 & \{1, 13, 89, 1325, 9077, 135137 ...\}\\
 194 & \{5, 13, 2897, 7561, 1686049, 4400489 ...\}\\
 433 & \{5, 29, 6466, 37666, 8399329, 48928105 ...\}\\
 169 & \{2, 29, 985, 14701, 499393, 7453378 ...\}\\
 89 & \{1, 34, 233, 9077, 62210, 2423525 ...\}\\
 1325 & \{13, 34, 51641, 135137, 205272962, 537169541 ...\}\\
 7561 & \{13, 194, 294685, 4400489, 6684339842, 99816291793 ...\}\\
\hline
\end{tabular}
\end{center}
\end{table}

 Thus, the Fr\"obenius Uniqueness Conjecture, which postulates that the region number of a \textit{Markov triplet} determines the remaining values of that triplet, is true. This statement follows from the analysis of the full Pell equation,
$$V[\{x,R,z\},n+1] ^2 - D(R) ~U[\{x,R,z\},n+1] ^2 = - (2R)^2,$$
 not the above, more dramatic,  \textit{Mathematica} analysis since the table could have been created by sorting the  $n = \{...~-3, -2, -1, 0, 1, 2~ ...\}$ values of 
 $U[\{x, R, z\}, n+1]$.

The Uniqueness Conjecture is fundamentally  an observation of a particular pattern within the Markov tree. The next section presents intricate, additional patterns derived from the triplet sequence functions.
 \newpage

 \section{Palindromic Repeat Cycles of Region Numbers along Edges}
\subsection{Triplet Parity Patterns}
Solutions to the Markov equation must be either all odd integers or two odd integers and one even integer. Thus, there are four possibilities for the parity of a triplet.
These four cases occur equally often within the Markov tree and at specific locations as shown in \mbox{Table 3.} Notice that although parity \textsc{types} 1, 2, 3, and 4 appear equally often on the left and right sides of the entire Markov tree, an individual region starting with a triplet of parity \textsc{type 4} produces only \textsc{type 3} children on its left edge and only \textsc{type 2} children on its right edge. The existence of repeating patterns shown here suggests that more interesting patterns may exist along left and right edges of a given region.

\begin{table}[ht]
\begin{center}
\caption{Parity Patterns of Left and Right Edge Members of Interior Region $R$}
\begin{tabular}{|c|c|c|}
\hline
Interior Region, $\{x, R, z\}$,					&Left Edge Members 		&Right Edge Members\\
Triplet Parity Pattern							&Repeating Parity Pattern		&Repeating Parity Pattern\\
and \textsc{Type}								&\textsc{Types}				&\textsc{Types}\\
\hline
\hline
1$\Longleftrightarrow$ \{\textsc{odd, odd, odd}\}	&\{\{4,2,1\}, \{4,2,1\}...\}	&\{\{4,3,1\}, \{4,3,1\}...\}\\
2$\Longleftrightarrow$ \{\textsc{even, odd, odd}\}	&\{\{2,1,4\}, \{2,1,4\}...\}	&\{\{1,4,3\}, \{1,4,3\}...\}\\
3$\Longleftrightarrow$ \{\textsc{odd, odd, even}\}	&\{\{1,4,2\}, \{1,4,2\}...\}	&\{\{3,1,4\}, \{3,1,4\}...\}\\
4$\Longleftrightarrow$ \{\textsc{odd, even, odd}\}	&\{\{3,3,3\}, \{3,3,3\}...\}	&\{\{2,2,2\}, \{2,2,2\}...\}\\
\hline
\end{tabular}
\end{center}
\end{table}

\subsection{Repeat Cycles of the Last Digits of Region Numbers}
The last one, two, three and four digits of the full set of Fibonacci numbers are known to have repeat cycle lengths of 60, 300, 1500, and 15000 terms. Table 4 shows the repeat cycle lengths of the last digits of the region numbers of triplets along the edges of the first few regions. The repeat cycle lengths of the odd-indexed Fibonacci and Pell numbers of \textsc{region 1} and \textsc{region 2} are half of the repeat cycle lengths of the full sets of Fibonacci and Pell numbers. 
Since the region numbers of the edges of a given region $\{x, R, z\}$ depend upon linear, constant coefficient, combinations of Lucas sequences $U_{k}(3R,1)$ and since these sequences have last digit repeat cycles, so do the region numbers.
Since the region numbers of the left edge of a region $\{x, R, z\}$ depend upon \mbox{$U[\{z, R, x\},n+1]$} while those of the right edge depend upon $U[\{x, R, z\},n+1]$ (which just has $x$ and $z$ interchanged), the left and right edges of a given interior region have the same repeat cycle lengths.

\begin{table}[ht]
\begin{center}
\caption{Repeat Cycles of Last Digits of Region Numbers Along Edges}
\begin{tabular}{|c|l|c|l|}
\hline
Markov			&Length of Repeat		&Markov 			&Length of Repeat\\
Region		 	&Cycles of Last			&Region			&Cycles of Last\\
Number			&\{1, 2, 3, 4\} digits		&Number			&\{1, 2, 3, 4\} digits\\
\hline
\hline
 1 & \{30, 150, 750, 7500\} & 2 &\{6, 30, 300, 3000\} \\
 5 &\{12, 60, 300, 1500\} & 13 &\{3, 15, 75, 750\} \\
 29 &\{15, 75, 375, 3750\} & 34 &\{5, 25, 500, 5000\} \\
 194 &\{5, 25, 500, 5000\} & 433 &\{3, 3, 30, 300\} \\
 169 &\{15, 75, 750, 7500\} & 89 &\{15, 75, 750, 7500\} \\
 1325 &\{12, 12, 60, 300\} & 7561 &\{30, 150, 750, 7500\} \\
 2897 &\{6, 30, 150, 1500\} & 6466 &\{10, 50, 500, 5000\} \\
 37666 &\{10, 50, 500, 5000\} & 14701 &\{30, 150, 750, 3750\} \\
 985 &\{12, 60, 300, 1500\} & 233 &\{3, 3, 30, 300\} \\
 9077 &\{6, 30, 150, 750\} & 135137 &\{6, 30, 150, 1500\} \\
\hline
\end{tabular}
\end{center}
\end{table}

There is no evidence that there are repeat cycles of the last digits of the Markov numbers, \{1, 2, 5, 13, 29 ...\}. However, some last digits occur more frequently than others. Fig.\ 2 shows the relative frequency of occurrence of the last single digit and last two digits for the first 2~million Markov numbers.
In the second graph of Fig.\ 2, only certain digit combinations appear. Even numbers with their last two digits of the form 2+20k and 18+20k occur least of all, followed by numbers with their last two digits of the form 10+20k, followed by numbers with their last two digits of the form 6+20k  and  14+20k. Odd numbers with their last two digits of the form 13+20k  and  17+20k occur next often, followed by numbers with their last two digits of the form 5+20k,  finally followed by numbers with their last two digits of the form 1+20k and  9+20k.\\ Note that the relative frequency of occurrence of even last digits is  \{3, 4, 5\}, while odd last digits show up triple that: \{9, 12, 15\}.

\begin{figure}[ht]
\begin{center}
\includegraphics[width=0.85\columnwidth]{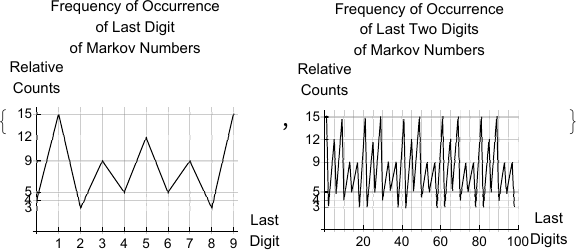}
\end{center}
\caption{Patterns in the frequency of occurrence of last digits of Markov numbers.}
\end{figure}

Table 4 showed that more than one triplet can have the same repeat cycle pattern and that a given pattern is somewhat correlated with the last digits of the triplet region number. For example, regions 16\underline{9} and 8\underline{9} share the same pattern. Table 5 shows the repeat cycle patterns for each of the possible two last digit patterns associated with Fig.\ 2. Note that the last two digits of the region number somewhat influence the repeat cycle patterns.

\begin{small}
\begin{table}[ht]
\begin{center}
\caption{Repeat Cycle Length Patterns vs. Region Number Last Two Digit Pattern}
\begin{footnotesize}
\begin{tabular}{|c|c|}
\hline
Pattern of Last		&Patterns of Length of Repeat Cycles\\
Two Digits of		&for Last \{1, 2, 3\} Digits\\
Region Number		&of Region Numbers\\
\hline
\hline
Mod[2 + 20k, 100]& \{6, 6, 12\}, \{6, 6, 60\}, \{6, 30, 300\}\\ 
Mod[18 + 20k, 100] & \{3, 3, 12\}, \{3, 3, 60\}, \{3, 15, 300\}\\
Mod[10 + 20k, 100] & \{4, 4, 4\}, \{4, 4, 20\}, \{4, 20, 100\}\\
Mod[6 + 20k, 100] & \{10, 50, 500\}\\ 
Mod[14 + 20k, 100] & \{5, 25, 500\}\\
Mod[13 + 20k, 100] & \{3, 3, 3\}, \{3, 3, 6\}, \{3, 3, 15\}, \{3, 3, 30\}, \{3, 15, 75\}\\
Mod[17 + 20k, 100] & \{6, 6, 6\}, \{6, 6, 30\}, \{6, 30, 150\}\\
Mod[5 + 20k, 100] & \{12, 12, 12\}, \{12, 12, 60\}, \{12, 60, 300\}\\
Mod[1 + 20k, 100] & \{30, 150, 750\}\\
Mod[9 + 20k, 100] & \{15, 75, 375\}, \{15, 75, 750\}\\
\hline
\end{tabular}
\end{footnotesize}
\end{center}
\end{table}
\end{small}
\subsection{Palindromic Repeat Cycles of Last Digits of Region Numbers}

Since the left and right edge region numbers have the same length repeat cycles  and since Table 1 showed that these repeat cycles are related, it is clear that appending the sequence of last digits of any region number repeat cycle of either side of any interior region to the corresponding sequence for the opposite side creates a \textit{palindromic} list.
Table 6 shows this palindromic behavior for the three digit repeat cycles of the first few interior regions listed in the first column. The second column shows the last three digits of the region numbers at the \{start ... end\} of the left edge three digit repeat cycle. The third column shows the corresponding last three digits of the region numbers at the \{start ... end\} of the right edge three digit repeat cycle. While repeat cycle lengths are independent of their starting point, palindromic repeat cycles are not. Region number palindromic sequences along edges start with $\textsc{H}_{sf} [\{x, R, z\},0].$ Thus, the first and last members of column two of this table are $x$ and $z$, then $z$ and $x$ for column three (all Mod 1000).

\begin{table}[ht]
\begin{center}
\caption{Palindromic Repeat Cycles of the Last Three Digits of Region Numbers}
\begin{footnotesize}
\begin{tabular}{|c|c|c|}
\hline
Region			&The Last Three Digits of	the Left	&The Last Three Digits of the Right\\	
$\{x, R, z\}$		&Edge Triplet Region Numbers		&Edge Triplet Region Numbers\\
				&\{start of edge ... end of edge\}	&\{start of edge ... end of edge\}\\
\hline
\hline
 \{1, 5, 2\} & \{1, 13, 194, 897...466, 433, 29, 2\} & \{2, 29, 433, 466...897, 194, 13, 1\}\\ 
\{1, 13, 5\} & \{1, 34, 325, 641...685, 561, 194, 5\} & \{5, 194, 561, 685...641, 325, 34, 1\}\\
     \{5, 29, 2\} & \{5, 433, 666, 509...818, 701, 169, 2\} & \{2, 169, 701, 818...509, 666, 433, 5\}\\ 
     \{1, 34, 13\} & \{1, 89, 77, 765...649, 137, 325, 13\} & \{13, 325, 137, 649...765, 77, 89, 1\}\\
     \{13, 194, 5\} & \{13, 561, 489, 37...621, 49, 897, 5\} & \{5, 897, 49, 621...37, 489, 561, 13\}\\ 
     \{5, 433, 29\} & \{5, 466, 329, 905...729, 105, 666, 29\} & \{29, 666, 105, 729...905, 329, 466, 5\}\\ 
     \{29, 169, 2\} & \{29, 701, 378, 945...266, 393, 985, 2\} & \{2, 985, 393, 266...945, 378, 701, 29\}\\
     \{1, 89, 34\} & \{1, 233, 210, 837...98, 525, 77, 34\} & \{34, 77, 525, 98...837, 210, 233, 1\}\\ 
     \{34, 1325, 13\} & \{34, 137, 541, 338...309, 962, 641, 13\} & \{13, 641, 962, 309...338, 541, 137, 34\}\\ 
     \{13, 7561, 194\} & \{13, 685, 842, 401...130, 793, 489, 194\} & \{194, 489, 793,130...401, 842, 685, 13\}\\ 
     \{194, 2897, 5\} & \{194, 49, 665, 466...825, 346, 261, 5\} & \{5, 261, 346, 825...466, 665, 49, 194\}\\
      \{5, 6466, 433\} & \{5, 557, 681, 481...253, 509, 329, 433\} & \{433, 329, 509, 253...481, 681, 557, 5\}\\ 
      \{433, 37666, 29\} & \{433, 105, 357, 181...585, 953, 509, 29\} & \{29, 509, 953, 585...181, 357, 105, 433\}\\ 
      \{29, 14701, 169\} & \{29, 818, 225, 357...417, 765, 378, 169\} & \{169, 378, 765,417...357, 225, 818, 29\}\\ 
  \hline
\end{tabular}
\end{footnotesize}
\end{center}
\end{table}

Interestingly, the palindromic sequences of these triplets have an \textit{internal structure.}
Consider the representative triplet $\{13, 7561, 194\}$ with repeat cycle pattern $\{30, 150, 750\}$ shown in Table 6. This repeat cycle pattern occurs more frequently than any other.
The two digit, 150-member repeat cycle (reading along the \textit{left} edge) has a 38-member \textit{palindromic} list, followed by a 112-member \textit{palindromic} list:\\

$\{\{13,85,42,1,41,2,25,73,34,49,33~...<16>...~33,49,34,73,25,2,41,1,42,85,13\},$\\
$\{94,89,93,30,97,21,46,97,5,18,89~...<90>...~89,18,5,97,46,21,97,30,93,89,94\}\}.$

The two digit, 150-member repeat cycle (reading along the \textit{right} edge) has the \textit{same} \mbox{112-member} \textit{palindromic} list, followed by the \textit{same} 38-member \textit{palindromic} list:\\
$\{\{94,89,93,30,97,21,46,97,5,18,89~...<90>...~89,18,5,97,46,21,97,30,93,89,94\},$\\
$\{\{13,85,42,1,41,2,25,73,34,49,33~...<16>...~33,49,34,73,25,2,41,1,42,85,13\}\}.$

Now consider another Table 6 representative triplet, $\{29, 14701, 169\},$ which has the same $\{30, 150, 750\}$ repeat cycle pattern. It has the following obvious property: Its last \textit{two} digit \textit{first palindromic} list along its left edge,
\begin{small}
$$\{29, 18, 25, 57, 46, 81, 97, 10, 33, 89, 34, 13, 5, 2, 1, 1, 2, 5, 13, 34, 89, 33, 10, 97, 81, 46, 57, 25, 18, 29\},$$
\end{small}
is composed entirely of the last two digits of the \textit{reversed and then forward-ordered first $15$ odd-indexed Fibonacci numbers (Mod100).} Looking back at the first representative triplet above, $\{13, 7561, 194\},$ it is clear that the digits of its \textit{first palindromic} list along its left edge, is composed entirely of the last single digits of \textit{reversed and then forward-ordered odd-indexed Fibonacci numbers.} 
Now consider the \textit{first palindromic} list along the left edge of the $1506^{th}$ non-singular triplet (the $483^{rd}$ member of depth 11 of the non-singular Markov tree) It has the same $\{30, 150, 750\}$ repeat cycle and its region number is: 
\begin{footnotesize} 
$$195307462239626228425885818346918767395108946306256881195001.$$
\end{footnotesize}
This list is composed of the last \textit{three} digits of the \textit{reversed and then forward-ordered first $47$ odd-indexed Fibonacci numbers (Mod1000)}:\\
\begin{small} 
$\{738, 309, 189, 258, 585, 497, 906, 221, 757, 50, 393, 129, 994, 853, 565, 842, 961, 41, 162, 445, 173, 74,\\
 49, 73, 170, 437, 141, 986, 817, 465, 578, 269, 229, 418, 25, 657, 946, 181, 597, 610, 233, 89, 34, 13, 5, 2, 1,\\
 1, 2, 5, 13, 34, 89, 233, 610, 597, 181, 946, 657, 25, 418, 229, 269, 578, 465, 817, 986, 141, 437, 170, 73, 49,\\
 74, 173, 445, 162, 41, 961, 842, 565, 853, 994, 129, 393, 50, 757, 221, 906, 497, 585, 258, 189, 309, 738\}.$\\
\end{small}

Although not shown, an analysis of the behavior of the \textit{second palindromic} lists along left edges of these triplets also shows Fibonacci numbers. 
 The single last digits of the second palindromic list of the triplet with region number 14701 is composed of the single last digits of \textit{ordered odd-indexed Fibonacci numbers.}  Similarly, the last two digits of the second palindromic list of the triplet above with the rather large region number $19530746...881195001$ is composed of the last two digits of \textit{ordered odd-indexed Fibonacci numbers.}

Now consider triplets with a different single edge repeat cycle pattern of $\{12, 60, 300\}$. This is the second most frequently occurring repeat cycle pattern. The palindromic sequences of \textit{all} of these triplets also have an \textit{internal structure.} The single digit lists of the sixty-member repeat cycle of each edge of these triplets has five twelve-member copies of the single digit list of the previous repeat cycle. These five twelve-member last digit lists are either ascending or descending cyclic \textit{Lucas numbers} (Mod 10), namely $\{1, 3, 4, 7, 1, 8, 9, 7, 6, 3, 9, 2\}.$  For example, the triplet \{1, 5, 2\} has five ascending copies of $\{1, 3, 4, 7, 1, 8, 9, 7, 6, 3, 9, 2\}$ along the left edge and five descending copies of $\{2, 9, 3, 6, 7, 9, 8, 1, 7, 4, 3, 1\}$ along the right edge. For these triplets, one edge is always ascending while the other is always descending because the first left edge member is always $x$ (Mod 10) and its last is always $z$ (Mod 10). The first right edge member is always $z$ (Mod 10) and its last is always $x$ (Mod 10).
Thus, only triplets whose $x$ and $z$ members' last single digits are adjacent members  of the cyclic Lucas number (Mod 10) list, $\{1, 3, 4, 7, 1, 8, 9, 7, 6, 3, 9, 2\},$ can have a $\{12, 60, 300\}$ repeat cycle. For example, $\{1, 5, 2\}$ works because $\{1, 2\}$ are adjacent cyclic Lucas numbers (Mod 10) and $\{169, 985, 2\}$ works because $\{9, 2\}$ are also adjacent. 

Other triplets which share repeat cycle patterns also have internal structure.

The Fibonacci and Pell numbers along the single edges of \textsc{Region 1} and \textsc{Region 2} also have palindromic repeat cycles. Tables 7 and 8 show the last few digits of the Fibonacci and Pell numbers at the end of different palindromic repeat cycles. Note that the leading digits shown in Table 7 settle down to the sequence of six numbers \{1,4,1,9,6,9\}. This pattern of six numbers repeats as the table is extended beyond the last shown entry of 89. Also note that \{1,4,1\} and \{9,6,9\} are tens complements of each other.

Interestingly, although the \textit{lengths} of repeat cycles of \textit{even-indexed} Fibonacci numbers are the same as the lengths of repeat cycles of \textit{odd-indexed} Fibonacci numbers, none of the \textit{even-indexed} Fibonacci number repeat cycles are palindromic since its first repeat cycle, 
$\{0, 1, 3, 8, 1, 5, 4, 7, 7, 4, 5, 1, 8, 3, 1, 0, 9, 7, 2, 9, 5, 6, 3,  3, 6, 5, 9, 2, 7, 9\},$
is not. Similarly, none of the \textit{even-indexed} Pell number repeat cycles are palindromic since its first,
$\{0, 2, 2, 0, 8, 8\},$
is also not. But note that these two repeat cycles have 10's complement patterns.

\begin{table}[ht]
\begin{center}
\caption{Last Digits of Region 1 Fibonacci Numbers at End of Repeat Cycle}
\begin{footnotesize}
\begin{tabular}{|c|r|r|r|r|r|r|r|}
\hline
Repeat Cycle&1&2&5&13&34&89\\
Length &&&&&&\\	
\hline
\hline
  30 & 41 & 62 & 45 & 73 & 74 & 49 \\
 150 & 201 & 802 & 205 & 813 & 234 & 889 \\
 750 & 1001 & 4002 & 1005 & 9013 & 6034 & 9089 \\
 7,500 & 10001 & 40002 & 10005 & 90013 & 60034 & 90089 \\
 75,000 & 100001 & 400002 & 100005 & 900013 & 600034 & 900089 \\
 750,000 & 1000001 & 4000002 & 1000005 & 9000013 & 6000034 & 9000089 \\
 7,500,000 & 10000001 & 40000002 & 10000005 & 90000013 & 60000034 & 90000089 \\
 75,000,000 & 100000001 & 400000002 & 100000005 & 900000013 & 600000034 & 900000089 \\
\hline
\end{tabular}
\end{footnotesize}
\end{center}
\end{table}
\begin{table}[ht]
\begin{center}
\caption{Last Digits of Region 2 Pell Numbers at End of Repeat Cycle}
\begin{footnotesize}
\begin{tabular}{|c|r|r|r|r|r|r|r|}
\hline
Repeat Cycle&1&5&29&169&985&5741\\
Length &&&&&&\\	
\hline
\hline
 6 & 41 & 85 & 69 & 29 & 5 & 1 \\
 30 & 701 & 905 & 729 & 469 & 85 & 41 \\
 300 & 7001 & 9005 & 7029 & 3169 & 1985 & 8741 \\
 3,000 & 70001 & 90005 & 70029 & 30169 & 10985 & 35741 \\
 30,000 & 700001 & 900005 & 700029 & 300169 & 100985 & 305741 \\
 300,000 & 7000001 & 9000005 & 7000029 & 3000169 & 1000985 & 3005741 \\
 3,000000 & 70000001 & 90000005 & 70000029 & 30000169 & 10000985 & 30005741 \\
 30,000000 & 700000001 & 900000005 & 700000029 & 300000169 & 100000985 & 300005741 \\
\hline
\end{tabular}
\end{footnotesize}
\end{center}
\end{table}

\pagebreak
\subsection{Farey Tree Sequence Functions}
The triplets of the Markov tree can be indexed by triplets of Farey numbers ranging from zero to one. Just as there are sequence functions for the region numbers along the left and right edges of all Markov regions, there are corresponding sequence functions for the Farey region numbers along the left and right edges of the Farey tree as shown in Table 9. The first two columns of the table give the Markov region and the corresponding Farey triplet at the head of that Markov region. The third column provides  the sequence functions of the Farey numbers corresponding to the Markov region numbers along the left and right edges of that region. Columns 4, 5 and 6 are analogous. Consider region 5. Its Farey triplet, $ \left\{0,\frac{1}{2},1\right\} $ corresponds to the Markov triplet $\{1,5,2\}$. The Markov region numbers along the left edge of this region, $\{13, 194, 2897...\}$, correspond to the Farey numbers $\frac{k}{2 k+1}$ for  $k=\{1,2,3...\}$, namely $\{\frac{1}{3},\frac{2}{5},\frac{3}{7}...\}.$ Analogously, the Markov region numbers along the right edge of this region, $\{29, 433, 6466...\}$, correspond to the Farey numbers $\frac{k+1}{2 k+1} $ for  $k=\{1,2,3...\}$, namely $\{\frac{2}{3},\frac{3}{5},\frac{4}{7}...\}.$
Note that these sequence functions have a simple pattern: If the Farey triplet at the head of a given region is denoted by 
$\{\frac{a}{b},\frac{x}{y},\frac{c}{d}\}$, then the \{left, right\} edge sequence functions of the Farey numbers corresponding to the Markov region numbers along the \{left, right\} edge of that region are given by
$\{\frac{x*k+a}{y*k+b},~ \frac{x*k+c}{y*k+d}\}.$ For large $k$, these sequence functions monotonically approach $\frac{x}{y}$ from below and above.

\begin{table}[ht]
\begin{center}
\caption{Sequence Functions of Farey Numbers Along Left \& Right Edge of a Region}
\begin{tabular}{|c|c|c|c|c|c|}
\hline
 Markov 			&Farey		&\{Left, Right\} 			&Markov		&Farey		&\{Left, Right\}\\
 Region			&Triplet at		&Edge Farey			&Region		&Triplet at		&Edge Farey\\
 Number			&Head of		& Region Number		&Number		&Head of		&Region Number\\
 				&Markov		&Sequence			&			&Markov		&Sequence\\
				&Region		&Function				&			&Region		&Function\\
				&			&k = \{1, 2, 3...\}		&			&			&k = \{1, 2, 3...\}\\
\hline
\hline
$ 5 $ & $ \{0,\frac{1}{2},1\} $ & $ \{\frac{k}{2 k+1},\frac{k+1}{2 k+1}\} $ & $ 13 $ & $ \{0,\frac{1}{3},\frac{1}{2}\} $ & $ \{\frac{k}{3 k+1},\frac{k+1}{3 k+2}\} $ \\
$ 29 $ & $ \{\frac{1}{2},\frac{2}{3},1\} $ & $ \{\frac{2 k+1}{3 k+2},\frac{2 k+1}{3 k+1}\} $ & $ 34 $ & $ \{0,\frac{1}{4},\frac{1}{3}\} $ & $ \{\frac{k}{4 k+1},\frac{k+1}{4 k+3}\} $ \\
$ 194 $ & $ \{\frac{1}{3},\frac{2}{5},\frac{1}{2}\} $ & $ \{\frac{2 k+1}{5 k+3},\frac{2 k+1}{5 k+2}\} $ & $ 433 $ & $ \{\frac{1}{2},\frac{3}{5},\frac{2}{3}\} $ & $ \{\frac{3 k+1}{5 k+2},\frac{3 k+2}{5 k+3}\} $ \\
$ 169 $ & $ \{\frac{2}{3},\frac{3}{4},1\} $ & $ \{\frac{3 k+2}{4 k+3},\frac{3 k+1}{4 k+1}\} $ & $ 89 $ & $ \{0,\frac{1}{5},\frac{1}{4}\} $ & $ \{\frac{k}{5 k+1},\frac{k+1}{5 k+4}\} $ \\
$ 1325 $ & $ \{\frac{1}{4},\frac{2}{7},\frac{1}{3}\} $ & $ \{\frac{2 k+1}{7 k+4},\frac{2 k+1}{7 k+3}\} $ & $ 7561 $ & $ \{\frac{1}{3},\frac{3}{8},\frac{2}{5}\} $ & $ \{\frac{3 k+1}{8 k+3},\frac{3 k+2}{8 k+5}\} $ \\
$ 2897 $ & $ \{\frac{2}{5},\frac{3}{7},\frac{1}{2}\} $ & $ \{\frac{3 k+2}{7 k+5},\frac{3 k+1}{7 k+2}\} $ & $ 6466 $ & $ \{\frac{1}{2},\frac{4}{7},\frac{3}{5}\} $ & $ \{\frac{4 k+1}{7 k+2},\frac{4 k+3}{7 k+5}\} $ \\
$ 37666 $ & $ \{\frac{3}{5},\frac{5}{8},\frac{2}{3}\} $ & $ \{\frac{5 k+3}{8 k+5},\frac{5 k+2}{8 k+3}\} $ & $ 14701 $ & $ \{\frac{2}{3},\frac{5}{7},\frac{3}{4}\} $ & $ \{\frac{5 k+2}{7 k+3},\frac{5 k+3}{7 k+4}\} $ \\
$ 985 $ & $ \{\frac{3}{4},\frac{4}{5},1\} $ & $ \{\frac{4 k+3}{5 k+4},\frac{4 k+1}{5 k+1}\} $ & $ 233 $ & $ \{0,\frac{1}{6},\frac{1}{5}\} $ & $ \{\frac{k}{6 k+1},\frac{k+1}{6 k+5}\} $ \\
\hline
\end{tabular}
\end{center}
\end{table}

\subsection{Markov Region Numbers vs. Associated Farey Numbers}

Figure 3 shows a plot of the Log[Markov Region Numbers] as a function of their associated Farey numbers. The plot is drawn as a series of lines connecting the $2^{d-1}$ points at a fixed depth, $d$, of the non-singular Markov tree. The plot includes lines for depth 2 through 8. The line at the bottom of the figure connects the two points, 
$\{\{\log[13], 1/3\},\{\log[29], 2/3\}\},$
 at depth two of the non-singular Markov tree, the next line connects the four points at depth three,
$\{\{\log[34], 1/4\},\{\log[194], 2/5\},\{\log[433], 3/5\},\{\log[169], 3/4\}\},$ and so on.
Since the region numbers within the Markov tree are not left-right symmetric,  this plot is not symmetric about the Farey number 1/2. 
Figure 4 shows a plot of the just the points of the Log[Markov Region Numbers] as a function of their associated Farey numbers for the first $2^{13}$ region numbers.

\begin{figure}[hbt!]
\begin{center}
\includegraphics[width=0.5\columnwidth]{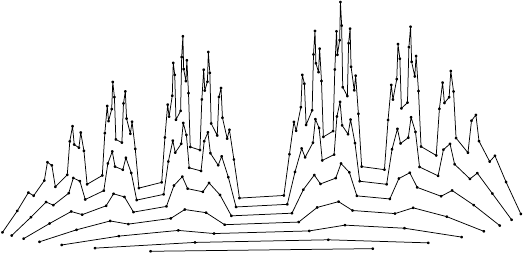}
\end{center}
\caption{Markov Region Numbers as a Function of their Farey Number}
\end{figure}
\begin{figure}[hbt!]
\begin{center}
\includegraphics[width=0.6\columnwidth]{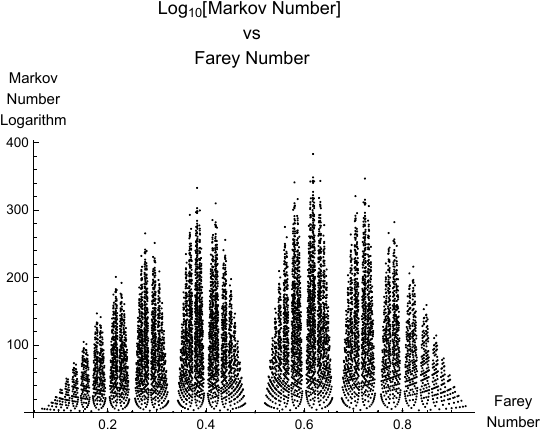}
\end{center}
\caption{Markov Region Numbers as a Function of their Farey Number}
\end{figure}

\section{All Markov Numbers are the Sum of Two Special Squares}
It is known that every Markov number is the sum of the squares of two integers. The wealth of intricate patterns amongst the members of a Markov Region~$R$ described above leads one to question what patterns might exist within these \textit{square terms}. Denote the smaller such integer square term by lower case $\sigma$ and the larger such integer square term by upper case $\Lambda$. Thus each Markov number, $R$, can be written as $\sigma^2 + \Lambda^2,$ almost invariably in multiple different ways. The following section analyzes cases of \textit{unique} square terms which sum to the Fibonacci region numbers of \textsc{Region~1}. These \textit{unique} cases have a signature pattern. This signature pattern provides the key to an algorithm which finds these special square terms for all triplets of all regions. The algorithm always finds \textit{one and only one} special pair of integers whose squares sum to the region number of interest. In order to find the special square terms for a region number of interest, it requires the complete triplet containing the region number of interest but also the special square terms of the grandparent of the triplet of interest, making the algorithm self-referential.

\subsection{The Signature Pattern of Special Square Terms of Region 1}
Since a given triplet, \(\{x, R, z\}\), can have multiple integer solutions to \mbox{$\sigma^2_{R}+\Lambda^2_{R}=R$,} this section analyzes only cases which have a \textit{unique} integer solution to $\sigma^2_{R}+\Lambda^2_{R}=R$ for several \textsc{Region 1} triplets as shown in Table 10. Any approach to understanding the integers which, when squared, sum to the Markov region numbers must encompass these cases.
\begin{table}[ht]
\begin{center}
\caption{\textsc{Region 1} Triplets which have  a \textit{Unique Solution} to $\sigma^2_{R}+\Lambda^2_{R}=R$}
\begin{tabular}{|c|c|c|c|c|c|c|}
\hline
\textsc{Region 1}	&Region						&Triplet		&Sibling						&Smallest											&Intermediate\\
Triplet:			&Square						&Sibling		&Square						&Triplet Member:									&Triplet Member:\\
\(\{x, R, z\}\)		&Terms:						&Number:		&Terms:						&$\pm (\sigma_{s} \Lambda_{R}-\Lambda_{s} \sigma_{R})$	&$\Lambda_{s} \Lambda_{R}+\sigma_{s} \sigma_{R}$\\
				&\{$\sigma_{R}$, $\Lambda_{R}$\}	&$3 x z - R$	&\{$\sigma_{s}$, $\Lambda_{s}$\}	&= Min$\{x, z\}$									& = Max$\{x, z\}$\\
\hline
\hline
\{1, 5, 2\} & \{1, 2\} & 1 & \{0, 1\} & $\sigma_{R} = 1$ & $\Lambda_{R}$ = 2 \\
 \{1, 13, 5\} & \{2, 3\} & 2 & \{1, 1\} & $\Lambda_{R} - \sigma_{R}$ = 1 & $\Lambda_{R} + \sigma_{R}$ = 5 \\
 \{1, 34, 13\} & \{3, 5\} & 5 & \{1, 2\} & $2 \sigma_{R} - \Lambda_{R}$ = 1 & $2 \Lambda_{R} + \sigma_{R}$ = 13 \\
 \{1, 89, 34\} & \{5, 8\} & 13 & \{2, 3\} & $2 \Lambda_{R} - 3 \sigma_{R}$ = 1 & $3 \Lambda_{R} + 2 \sigma_{R}$ = 34 \\
 \{1, 233, 89\} & \{8, 13\} & 34 & \{3, 5\} & $5 \sigma_{R} - 3 \Lambda_{R}$ = 1 & $5 \Lambda_{R} + 3 \sigma_{R}$ = 89 \\
\hline
\end{tabular}
\end{center}
\end{table}

Column 1 of this table lists the \textsc{Region 1} triplet of interest. Column 2 gives the  \textit{unique} integers whose squares sum to that triplet's region number. Column 3 gives the secondary solution to the Markov equation for fixed $x$ and $z$, namely $s = 3 x z -R,$ the region sibling number.
Column 4 gives the  \textit{unique} integers whose squares sum to that sibling number. 
The last two columns of this table are key: They give the simple multiples of $\sigma_{R}$ and $\Lambda_{R}$ which form the smallest and intermediate triplet members: Min$\{x,z\}$ and Max$\{x,z\}$. Aside from signs, these multiples are the unique integers whose squares sum to the region sibling number, \{$\sigma_{s}$, $\Lambda_{s}$\}.
These terms are $\pm (\sigma_{s} \Lambda_{R}-\Lambda_{s} \sigma_{R})$ for the smallest triplet member, Min$\{x,z\}$, and ($\Lambda_{s} \Lambda_{R}+\sigma_{s} \sigma_{R}$) for the intermediate triplet member, Max$\{x,z\}$.

These terms are an instance of the Brahmagupta-Fibonacci identity since
$$R s =  x^2 + z^2$$
can be rewritten as:
$$( \Lambda_{R}^2 + \sigma_{R}^2)( \Lambda_{s}^2 + \sigma_{s}^2) = 
  (\Lambda_{s} \Lambda_{R}+\sigma_{s} \sigma_{R})^2 + (\sigma_{s} \Lambda_{R}-\Lambda_{s} \sigma_{R})^2.$$

The results in these last two columns suggest an algorithm for finding the unique special square terms \{$\sigma_{R}$, $\Lambda_{R}$\} which, when squared,  sum to the region number, $R$, of any triplet, \(\{x, R, z\}\), avoiding multiple solutions which typically arise from solving the integer domain equation: $\sigma^2_{R}+\Lambda^2_{R} = R$. 
When this algorithm is applied to the edges of an arbitrary region, it is found that these special square terms have simple  recurrence kernels, generating functions and sequence functions.
This algorithm is described next.

 \subsection{An Algorithm for Evaluating the Special Square Terms of Any Triplet} 
 The algorithm, named $Q [\{x, R, z\}]$, requires a list of the ordered, singular and non-singular Markov triplets:
 
  \textsc{MarkovList} = $\{\{1,1,1\},\{1,2,1\},\{1,5,2\},\{1,13,5\},\{5,29,2\},\{1,34,13\}...\}.$
  
  It accepts as input any triplet, \(\{x, R, z\}\), within \textsc{MarkovList} and outputs the special square terms of both the triplet's region number, $R$, and its region sibling number, $s$. It accomplishes this by first computing the sibling number, $s = 3 x z-R,$ corresponding to the input triplet and retrieving the particular triplet which has a region number of $s$ from \textsc{MarkovList}. The algorithm then recursively calls itself to evaluate the special square terms, $\{\sigma_{s},  \Lambda_{s}\}$, of this triplet with region number, $s < R$. It then forms and solves the two linear equations contained in the last two columns of Table 10 for the special square terms of the input region number, $\{ \sigma_{R},  \Lambda_{R}\}$:

$$\pm (\sigma_{s} \Lambda_{R} - \Lambda_{s} \sigma_{R}) = \textsc{min} \{x, z\}$$
$$\Lambda_{s} \Lambda_{R} + \sigma_{s} \sigma_{R} = \textsc{max} \{x, z\}$$
(once the appropriate $\pm$ sign is determined). Finally, the algorithm outputs the special square terms for both the region number and the region sibling number. 
Since the algorithm is self-referential, it requires three startup values to avoid an infinite loop.

These startup values are what would be the ordered output, $\{\{ \sigma_{R},  \Lambda_{R}\}, \{ \sigma_{s},  \Lambda_{s}\}\},$ of the algorithm for the first three  \textsc{MarkovList} triplets:
  $$Q [\{1, 1, 1\}] =\{\{0,1\},\{1,1\}\}$$
  $$Q [\{1, 2, 1\}] =\{\{1,1\},\{0,1\}\}$$
  $$Q [\{1, 5, 2\}] =\{\{1,2\},\{0,1\}\}.$$

In detail, the algorithm proceeds as follows:
\begin{itemize}
\item{Accept Input: $\{x, R, z\}$ (\textit{The Triplet of Interest})}
\item{Calculate: $s = (3 x z-R)$ (\textit{The Input Triplet's Region Sibling Number})}
\item{Retrieve: $\{x', s, z'\}$ from \textsc{MarkovList} (\textit{The Triplet with Region Number, $s$})}
\item{Calculate: $\{\sigma_{s},  \Lambda_{s}\} = \textsc{part 1 of } Q[\{x', s, z'\}]$ (\textit{Special Squares of Sibling Number})}
\item{Calculate: \textsc{Sign} = \textsc{RegionSign} [$\{x, R, z\}$] (\textit{The $\pm$ sign for the} \textsc{SmallestTerm})}
\item{Calculate: \textsc{SmallestTerm} = \textsc{Sign} $*$ ($\sigma_{s} \Lambda_{R} - \Lambda_{s} \sigma_{R}$) (\textit{The Smallest Triplet Member})}
\item{Calculate: \textsc{IntermediateTerm} = $\Lambda_{s} \Lambda_{R} + \sigma_{s} \sigma_{R}$ (\textit{The Intermediate Triplet Member})}
\item{Solve: \{\textsc{SmallestTerm} = \textsc{Min}[\{x, z\}], \textsc{IntermediateTerm} = \textsc{Max}[\{x, z\}]\} for the variables: $\{\sigma_{R}, \Lambda_{R}\}$ (\textit{The two Special Square Terms whose squares sum to the Triplet Region Number}, $R$.}
\item{Return Output: $\{\{\sigma_{R},  \Lambda_{R}\}, \{\sigma_{s},  \Lambda_{s}\}\}$ (\textit{The Region and Sibling Number Square terms})}
\end{itemize}

With two exceptions required for the singular Markov regions, the function \textsc{RegionSign} is defined as the product of  two functions \textsc{LeftRightSign} and \textsc{ParitySign}:

$\textsc{RegionSign} [\{1, 2, 1\}] =-1$\\
$\textsc{RegionSign} [\{1, 5, 2\}] =-1$\\
$\textsc{RegionSign} [\{x, R, z\}] = \textsc{LeftRightSign} [\{x, R, z\}] *\textsc{ParitySign} [\{x, R, z\}]$

$\textsc{LeftRightSign} [\{x, R, z\}] = +1$ if the position of the triplet $\{x, R, z\}$ within \textsc{MarkovList} is a multiple of two and $-1$ otherwise.

$\textsc{ParitySign} [\{x, R, z\}] = -1$ if the position of the triplet $\{x, R, z\}$ within \textsc{MarkovList} is a multiple of three and $+1$ otherwise.

Table 11 gives the special square terms output from this algorithm for triplets along the left and right edges of sample \textsc{Region 5.} 
\begin{table}[ht]
\caption{Special Square Terms Along Left And Right Edges of Sample \textsc{Region 5}}
\begin{small}
\begin{tabular}{|c|c|c|c|}
\hline
Left Edge			&Special Square				&Right Edge 			&Special Square\\
Region 5 Triplet 	&Terms of the					&Region 5 Triplet		&Terms of the\\
$\{x, R, z\}$		&Left EdgeTriplet's				&$\{x, R, z\}$			&Right EdgeTriplet's\\
				&Region Number				&					&Region Number\\
				&$\{\sigma_{R}, \Lambda_{R}\}\rightarrow\sigma_{R}^2+\Lambda_{R}^2$	&					
				&$\{\sigma_{R}, \Lambda_{R}\}\rightarrow\sigma_{R}^2+\Lambda_{R}^2$\\
\hline
\hline
\{1,\;13,\;5\} & \{2,\;3\} $\rightarrow$ 13 & \{5,\;29,\;2\} & \{2,\;5\} $\rightarrow$ 29\\
 \{13,\;194,\;5\} & \{5,\;13\} $\rightarrow$ 194 & \{5,\;433,\;29\} & \{12,\;17\} $\rightarrow$ 433\\
 \{194,\;2897,\;5\} & \{31,\;44\} $\rightarrow$ 2897 & \{5,\;6466,\;433\} & \{29,\;75\} $\rightarrow$ 6466\\
 \{2897,\;43261,\;5\} & \{75,\;194\} $\rightarrow$ 43261 & \{5,\;96557,\;6466\} & \{179,\;254\} $\rightarrow$ 96557\\
 \{43261,\;646018,\;5\} & \{463,\;657\} $\rightarrow$ 646018 & \{5,\;1441889,\;96557\} & \{433,\;1120\} $\rightarrow$ 1441889\\
\hline
\end{tabular}
\end{small}
\end{table}

By analyzing the output of this algorithm along left and right edges of different regions, it is found that
the \{odd,  even\} indexed square terms along an edge naturally group together. Each group has a recurrence kernel of $\{3R,~-1\}$. Thus, generating functions and sequence functions for the special square terms along an edge of a region $\{x, R, z\}$  are given in terms of four two-member lists $\{\alpha, \beta, \gamma, \delta\}$ by:
$$
\textsc{K}_{gf} [\{x, R, z\},n] = 
\left(
\begin{array}{cc}
 \frac{\alpha - \gamma \; n} {1-3 R n +n^2} & \textsc{for odd } n \\
 ~\\
 \frac{\beta - \delta \; n} {1-3 R n +n^2} & \textsc{for even } n \\
\end{array}
\right.
$$
$$
\textsc{K}_{sf} [\{x, R, z\},n] = 
\left(
\begin{array}{cc}
 \alpha \; U_{n}(3R,1)-\gamma \; U_{n-1}(3R,1) & \textsc{for odd } n \\
 \beta \; U_{n}(3R,1)-\delta \; U_{n-1}(3R,1) & \textsc{ for even } n \\
\end{array}
\right.
$$
The functions  $\textsc{K}_{gf} [\{x, R, z\},n]$ and $\textsc{K}_{sf} [\{x, R, z\},n]$ give the ordered list $\{\sigma_{R},  \Lambda_{R}\}$ for $n \ge 1$.
The four lists $\{\alpha, \beta, \gamma, \delta\}$ are computed as follows:
\newpage
\begin{itemize}
\item{Accept Input: $\{triplet_{1}, triplet_{2}\}$ (\textit{The First Two Triplets along the Edge of Interest})}
\item{Evaluate: (\textit{Special Square Terms of these two input triplets})\\
$$\{\{\{\sigma_{R},  \Lambda_{R}\}_{1},\{\sigma_{s},  \Lambda_{s}\}_{1}\},\{\{\sigma_{R},  \Lambda_{R}\}_{2},\{\sigma_{s},  \Lambda_{s}\}_{2}\} \}\} = \{Q[triplet_{1}],Q[triplet_{2}]\}.$$ }
\item{Evaluate: \textit{The function} \textsc{RegionSign} \textit{using the first input triplet}
$$sign_{1} =\textsc{RegionSign}[triplet_{1}].$$}
\item{Form the lists $\{\alpha, \beta\}:$ (\textit{The special square terms of the region numbers})
$$\{\alpha, \beta\} = \{\{\sigma_{R},  \Lambda_{R}\}_{1},\{\sigma_{R},  \Lambda_{R}\}_{2}\}.$$}
\item{Form the lists $\{\gamma, \delta\}:$ (\textit{The special square terms of the sibling numbers})
$$\{\gamma, \delta\} = \{sign_1 *\{-\Lambda_{s},  \sigma_{s}\}_{1},\{\sigma_{s},  \Lambda_{s}\}_{2}\}.$$}
\item{Output the lists: $\{\alpha, \beta,\gamma,\delta\}.$}
\end{itemize}

Tables 12 and 13 give the left and right edge lists $\{\alpha, \beta, \gamma, \delta\}$ for a few Markov regions.
\begin{table}[ht]
\begin{center}
\caption{\textit{Left Edge} Lists $\{\alpha, \beta,\gamma,\delta\}$ for Region $\{x, R, z\}$}
\begin{footnotesize}
\begin{tabular}{|c|c|c|c|c|c|c|}
\hline
\textsc{Region} \(\{x, R, z\}\)	&List	 $\alpha$		&List	 $\beta$		&List	 $\gamma$		&List $\delta$	\\
\hline
\{1, 5, 2\}		& \{2, 3\}		& \{5, 13\}		& \{-1, 1\}		& \{0, 1\}\\ 
\{1, 13, 5\}		& \{3, 5\}		& \{13, 34\}	& \{2, -1\}		& \{0, 1\}\\ 
\{5, 29, 2\}		& \{12, 17\}	& \{75, 179\}	& \{-1, 1\}		& \{1, 2\}\\ 
\{1, 34, 13\}	& \{5, 8\}		& \{34, 89\}	& \{-3, 2\}		& \{0, 1\}\\ 
\{13, 194, 5\}	& \{44, 75\}	& \{1208, 1715\}& \{2, -1\}		& \{2, 3\}\\ 
\{5, 433, 29\}	& \{29, 75\}	& \{1120, 2673\}& \{-5, 2\}		& \{1, 2\}\\ 
\{29, 169, 2\}	& \{70, 99\}	& \{1043, 2523\}& \{-1, 1\}		& \{2, 5\}\\
\hline
\end{tabular}
\end{footnotesize}
\end{center}
\end{table}
\begin{table}[ht]
\begin{center}
\caption{\textit{Right Edge} Lists $\{\alpha, \beta,\gamma,\delta\}$ for Region $\{x, R, z\}$}
\begin{footnotesize}
\begin{tabular}{|c|c|c|c|c|c|c|}
\hline
\textsc{Region} \(\{x, R, z\}\)	&List	 $\alpha$		&List	 $\beta$		&List	 $\gamma$		&List $\delta$	\\
\hline
\{1, 5, 2\}		& \{2, 5\}		& \{12, 17\}		& \{1, 0\}	& \{1, 1\}\\ 
\{1, 13, 5\}		& \{5, 13\}		& \{44, 75\}		& \{1, 0\}	& \{1, 2\}\\
\{5, 29, 2\}		& \{5, 12\}		& \{70, 99\}		& \{-2, 1\}	& \{1, 1\}\\ 
\{1, 34, 13\}	& \{13, 34\}	& \{196, 311\}		& \{1, 0\}	& \{2, 3\}\\
\{13, 194, 5\}	& \{31, 44\}	& \{657, 1120\}		& \{3, -2\}	& \{1, 2\}\\ 
\{5, 433, 29\}	& \{75, 179\}	& \{2523, 6524\}	& \{-2, 1\}	& \{2, 5\}\\ 
\{29, 169, 2\}	& \{12, 29\}	& \{408, 577\}		& \{5, -2\}	& \{1, 1\}\\
\hline
\end{tabular}
\end{footnotesize}
\end{center}
\end{table}

\begin{itemize}
\item{$\textsc{K}_{sf} [\{x, R, z\},n]$ \textbf{Works for All Regions and All Integer $n$}

These sequence functions for the special square terms were derived for positive values of $n$.
However, the \{left, right\} edge expression for odd negative $n$ produces values for the \{right, left\} edge expression for even positive $n$. Similarly, the \{left, right\} edge expression for even negative $n$ produces values for \{right, left\} edge expression for odd positive $n$. The smaller square terms are swapped with the larger terms and consistent negative signs appear, but $\sigma_R^2 + \Lambda_R^2$ is preserved.  Each edge expression has the full information content of all triplets along both edges although only two triplets along one edge were sufficient to derive the sequence functions. This behavior occurs for all interior regions. Table 14 shows these results for sample \textsc{Region 5}.
\begin{table}[ht]
\begin{center}
\caption{Left and Right Edge Outputs for Positive and Negative Arguments to the Sequence Function, $\textsc{K}_{sf} [\{x, R, z\},n]$, for the Special Square Terms  of  Sample \textsc{Region 5}}
\begin{tabular}{|c|c|c||c|c|c|c|}
\hline
\textsc{Index}	&Left Edge			&Right Edge			&Left Edge			&Right Edge\\
$|n|$			&Expression			&Expression			&Expression			&Expression\\
			&Output for			&Output for			&Output for			&Output for\\
			&\textsc{positive} $n$	&\textsc{positive} $n$	&\textsc{negative} $n$	&\textsc{negative} $n$\\
\hline
\hline
Odd&&&&\\
1	&\{2,\;3\} & \{2,\;5\} & \{-17,\;12\} & \{13,\;-5\} \\
 3	&\{31,\;44\} & \{29,\;75\} & \{-254,\;179\} & \{194,\;-75\} \\
 5	&\{463,\;657\} & \{433,\;1120\} & \{-3793,\;2673\} & \{2897,\;-1120\} \\
 \hline
 \hline
 Even&&&&\\
 2	&\{5,\;13\} & \{12,\;17\} & \{-5,\;2\} & \{3,\;-2\} \\
 4	&\{75,\;194\} & \{179,\;254\} & \{-75,\;29\} & \{44,\;-31\} \\
 6	&\{1120,\;2897\} & \{2673,\;3793\} & \{-1120,\;433\} & \{657,\;-463\} \\
\hline
\end{tabular}
\end{center}
\end{table}
}
\end{itemize}

\section{Palindromic Repeat Cycles of Special Square Terms}

\subsection{Repeat Cycles of Last Digits of Special Square Terms}
The last one, two, three and four digits of the region numbers of triplets along the left and right edges of Markov regions were shown to have repeat cycles in Section 3.2. Here, the situation is more complicated. The special square terms consist of a \textit{pair} of numbers which naturally form \textit{two} groups with odd and even indices. 
Table 15 shows the repeat cycle lengths for the square terms of the single edge Fibonacci and Pell regions followed by the repeat cycle lengths of the square terms of interior Markov regions.  
\begin{table}[ht]
\begin{center}
\caption{Repeat Cycle Patterns of Last Digits of Square Terms for Markov Regions}
\begin{tabular}{|c|l|c|l|}
\hline
Markov			&Length of Repeat		&Markov 			&Length of Repeat\\
Region		 	&Cycles of Last			&Region			&Cycles of Last\\
Number			&\{1, 2, 3, 4\} digits		&Number			&\{1, 2, 3, 4\} digits\\
\hline
\hline
 1 & \{30, 150, 750, 7500\} & 2 &\{6, 30, 300, 3000\} \\
 5 &\{12, 60, 300, 1500\} & 13 &\{3, 15, 75, 750\} \\
 29 &\{15, 75, 375, 3750\} & 34 &\{5, 25, 500, 5000\} \\
 194 &\{5, 25, 500, 5000\} & 433 &\{3, 3, 30, 300\} \\
 169 &\{15, 75, 750, 7500\} & 89 &\{15, 75, 750, 7500\} \\
 1325 &\{12, 12, 60, 300\} & 7561 &\{30, 150, 750, 7500\} \\
 2897 &\{6, 30, 150, 1500\} & 6466 &\{10, 50, 500, 5000\} \\
 37666 &\{10, 50, 500, 5000\} & 14701 &\{30, 150, 750, 3750\} \\
 985 &\{12, 60, 300, 1500\} & 233 &\{3, 3, 30, 300\} \\
 9077 &\{6, 30, 150, 750\} & 135137 &\{6, 30, 150, 1500\} \\
\hline
\end{tabular}
\end{center}
\end{table}

The repeat cycle length patterns shown here apply for either choice of left edge or right edge, either choice of odd-indexed or even-indexed and either choice of $\sigma_R$ or $\Lambda_R$. All share the same recurrence kernel, $\{3R, -1\},$ and thus the same repeat cycle patterns. In addition, the repeat cycle patterns shown here are identical to the repeat cycle patterns of the region numbers shown in Table 4 (which also have the recurrence kernel $\{3R, -1\}$).

\subsection{Palindromic Repeat Cycles of Last Digits of Special Square Terms}
From Table 14, it is clear that joining the sequence of last digits of any repeat cycle of region number special square terms, $\{\sigma_R,\Lambda_R\},$  of either side of any interior region to the corresponding sequence for the opposite side creates a somewhat unusual palindromic list.
Table 16 shows this behavior for the two digit repeat cycles of the first few interior regions. It was created by shortening the three digit repeat cycles. 
The palindromic repeat cycles of the region number, $R$, described in Section 3.3, started with $\textsc{H}_{sf} [\{x, R, z\},0].$ Here, the palindromic repeat cycles of the region number special square terms, $\{\sigma_R, \Lambda_R\}$, must start at the same place. 

\begin{table}[ht]
\small
\begin{center}
\caption{Left and Right Edge Palindromic Repeat Cycles of Last Digits of Special Squares}
\begin{footnotesize}
\begin{tabular}{|c||c|c|}
\hline
 Region 		&The Last Two Digits of the First			&The Last Two Digits of the Last\\	
$\{x,R,z\}$		&Few Left Edge Special Squares,			&Few Right Edge Special Squares\\	
			& $\{\sigma_R, \Lambda_R\},$ At the Start	&$\{\sigma_R, \Lambda_R\},$ At the End\\	
			&of the Repeat Cycle					&of the Repeat Cycle\\	
\hline
\hline
    \{1,5,2\}&\{0,\underline{1}\},\{2,\underline{3}\},\{5,\underline{13}\},\{31,\underline{44}\},\{75,\underline{94}\}..&..\{\underline{94},25\},\{\underline{44},69\},\{\underline{13},95\},\{\underline{3},98\},\{\underline{1},0\}\\
\{1,13,5\}&\{0,\underline{1}\},\{\underline{3},5\},\{13,\underline{34}\},\{\underline{15},96\},\{7,\underline{25}\}..&..\{\underline{25},93\},\{4,\underline{15}\},\{\underline{34},87\},\{95,\underline{3}\},\{\underline{1},0\}\\
\{5,29,2\}&\{\underline{1},2\},\{12,\underline{17}\},\{\underline{75},79\},\{45,\underline{78}\},\{\underline{24},71\}..&..\{29,\underline{24}\},\{\underline{78},55\},\{21,\underline{75}\},\{\underline{17},88\},\{98,\underline{1}\}\\
  \{1,34,13\}&\{0,\underline{1}\},\{5,\underline{8}\},\{34,\underline{89}\},\{13,\underline{14}\},\{68,\underline{77}\}..&..\{\underline{77},32\},\{\underline{14},87\},\{\underline{89},66\},\{\underline{8},95\},\{\underline{1},0\}\\
\{13,194,5\}&\{2,\underline{3}\},\{\underline{44},75\},\{8,\underline{15}\},\{\underline{6},51\},\{54,\underline{27}\}..&..\{\underline{27},46\},\{49,\underline{6}\},\{\underline{15},92\},\{25,\underline{44}\},\{\underline{3},98\}\\
\{5,433,29\}&\{\underline{1},2\},\{29,\underline{75}\},\{\underline{20},73\},\{76,\underline{23}\},\{\underline{79},25\}..&..\{75,\underline{79}\},\{\underline{23},24\},\{27,\underline{20}\},\{\underline{75},71\},\{98,\underline{1}\}\\
\{29,169,2\}&\{2,\underline{5}\},\{70,\underline{99}\},\{43,\underline{23}\},\{91,\underline{92}\},\{99,\underline{56}\}..&..\{\underline{56},1\},\{\underline{92},9\},\{\underline{23},57\},\{\underline{99},30\},\{\underline{5},98\}\\
\{1,89,34\}&\{0,\underline{1}\},\{\underline{8},13\},\{89,\underline{33}\},\{\underline{31},74\},\{63,\underline{10}\}..&..\{\underline{10},37\},\{26,\underline{31}\},\{\underline{33},11\},\{87,\underline{8}\},\{\underline{1},0\}\\
\{34,1325,13\}&\{\underline{3},5\},\{96,\underline{11}\},\{\underline{29},90\},\{3,\underline{23}\},\{\underline{72},45\}..&..\{55,\underline{72}\},\{\underline{23},97\},\{10,\underline{29}\},\{\underline{11},4\},\{95,\underline{3}\}\\
\{13,7561,194\}&\{2,\underline{3}\},\{94,\underline{7}\},\{81,\underline{41}\},\{15,\underline{76}\},\{21,\underline{0}\}..&..\{\underline{0},79\},\{\underline{76},85\},\{\underline{41},19\},\{\underline{7},6\},\{\underline{3},98\}\\
\hline
\end{tabular}
\end{footnotesize}
\end{center}
\end{table}

The palindromic behavior within this table is unusual since it mimics  the behavior of the sequence functions for the special square terms shown in Table 14. There, a single sequence function was evaluated for positive and negative arguments along both edges of a region. That evaluation systematically \textit{interchanged} the $\{\sigma_R, \Lambda_R\}$ terms and introduced \textit{negative signs.} If the non-underlined terms in Table 16 are ignored, the remaining terms have a proper palindromic behavior. However, this \textit{interchanges} the $\{\sigma_R, \Lambda_R\}$ terms (sometimes in the odd-indexed \textit{and} even-indexed terms, other times in the odd-indexed \textit{or} even-indexed terms). Alternately, if the underlined terms in this table are ignored and the 100's complement of just one side of the table is formed, then the remaining terms also have a proper palindromic behavior. Note that this \textit{interchanges} the $\{\sigma_R, \Lambda_R\}$ terms as before.

The special squares of the Fibonacci numbers in \textsc{Region 1} also exhibit this unusual behavior within their single edge repeat cycles. Not surprisingly, the odd-indexed $\sigma_R$ and even-indexed $\Lambda_R$ terms of \textsc{Region 1} are odd-indexed Fibonacci numbers. But the even-indexed $\sigma_R$ and odd-indexed $\Lambda_R$ terms of \textsc{Region 1} are \textit{even-indexed} Fibonacci numbers.

Section 3.3 analyzed the \textit{internal structure} of the palindromic repeat cycles of last digits of region numbers. There, it was found that the lists of last digits of region numbers of the region $\{29, 14701, 169\}$ was composed of patterns of odd-indexed Fibonacci numbers. Here, the lists of palindromic repeat cycles of last digits of special square terms display similar behavior for the same $\{29, 14701, 169\}$ region. Since the special square terms are a pair of numbers which naturally fall into odd- and even-indexed groups, the behavior is more complicated. 
The lists of single last digits of the 30-member repeat cycles of both odd- and even-indexed $\Lambda_R$ terms are composed of two different cyclic rotations of the first 30 \textit{odd-indexed} Fibonacci numbers (Mod 10). 
The corresponding lists of single last digits of the 30-member repeat cycles of both odd- and even-indexed $\sigma_R$ terms are composed of two different cyclic rotations of the first 30 \textit{even-indexed} Fibonacci numbers (Mod 10).

\subsection{Oscillatory Behavior of Special Square Terms}
Fig.\ 5 displays the surprising oscillatory behavior of the special square terms along each edge of a typical region, here \textsc{region 5}. This behavior does not diminish along an edge and occurs for all edges of all regions. In each instance, the ratio of the upper to lower oscillation bound analytically converges to a limit determined by the particular region. 
For the Fibonacci region, the ratio of the upper oscillation bound to the lower oscillation bound analytically converges to the golden ratio, $(1+\sqrt{5})/2$, while for the Pell region this ratio analytically converges to the silver ratio, $(1+\sqrt{2}).$
\begin{figure}[ht]
\begin{center}
\includegraphics[width=0.6\columnwidth]{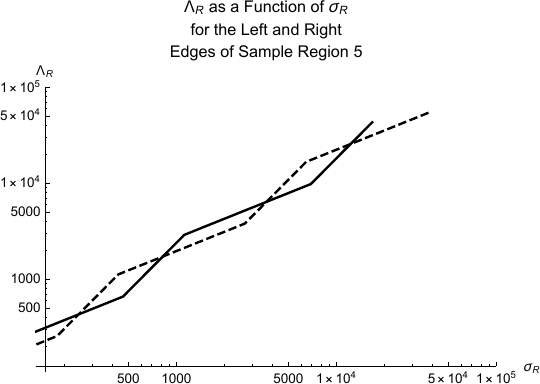}
\end{center}
\caption{The Oscillatory Behavior of $\Lambda_R$  as a Function of $\sigma_R$ along the Left and Right Edges of Sample \textsc{Region 5}.}
\end{figure}

\section{Inspiration}
This work was begun after reading an article in Quanta Magazine. The abundance of patterns of the Fibonacci and Pell numbers in \textsc{Region 1} and \textsc{Region 2} of the Markov tree provided signposts for further analysis using \textit{Mathematica}. 

The hope is that others will find these results as interesting as I have.

\section{References}
Kramer, A., \& Kramer, A. (2023, December 12). A triplet tree forms one of the most beautiful structures in math. Quanta Magazine. https://www.quantamagazine.org/a-triplet-tree-forms-one-of-the-most-beautiful-structures-in-math-20231212/

\end{document}